
\input amstex.tex
\documentstyle{amsppt}
\magnification1200
\hsize=12.5cm
\vsize=18cm
\hoffset=1cm
\voffset=2cm

\def\DJ{\leavevmode\setbox0=\hbox{D}\kern0pt\rlap
{\kern.04em\raise.188\ht0\hbox{-}}D}
\def\dj{\leavevmode
 \setbox0=\hbox{d}\kern0pt\rlap{\kern.215em\raise.46\ht0\hbox{-}}d}

\def\txt#1{{\textstyle{#1}}}
\baselineskip=13pt
\def\hf{{\textstyle{1\over2}}}
\def\a{\alpha}\def\b{\beta}
\def\d{{\,\roman d}}
\def\e{\varepsilon}
\def\D{\Delta}
\def\b{\beta}
\def\no{\noindent}

\def\G{\Gamma}

\def\s{\sigma}
\def\t{\theta}
\def\={\;=\;}

\def\zt{\zeta(\hf+it)}

\def\D{\Delta}
\def\no{\noindent}
\def\R{\Re{\roman e}\,} 
\def\z{\zeta}

\def\t{\theta}
\def\hf{{\textstyle{1\over2}}}
\def\txt#1{{\textstyle{#1}}}

\def\le{\leqslant}
\def\ge{\geqslant}
\font\tenmsb=msbm10
\font\sevenmsb=msbm7
\font\fivemsb=msbm5
\newfam\msbfam
\textfont\msbfam=\tenmsb
\scriptfont\msbfam=\sevenmsb
\scriptscriptfont\msbfam=\fivemsb
\def\Bbb#1{{\fam\msbfam #1}}

\def \NN {\Bbb N}

\def \RR {\Bbb R}
\def \ZZ {\Bbb Z}

\font\ff=cmr8
\def\txt#1{{\textstyle{#1}}}
\baselineskip=13pt

\font\teneufm=eufm10
\font\seveneufm=eufm7
\font\fiveeufm=eufm5
\newfam\eufmfam
\textfont\eufmfam=\teneufm
\scriptfont\eufmfam=\seveneufm
\scriptscriptfont\eufmfam=\fiveeufm
\def\mathfrak#1{{\fam\eufmfam\relax#1}}

\font\tenmsb=msbm10
\font\sevenmsb=msbm7
\font\fivemsb=msbm5
\newfam\msbfam
     \textfont\msbfam=\tenmsb
      \scriptfont\msbfam=\sevenmsb
      \scriptscriptfont\msbfam=\fivemsb
\def\Bbb#1{{\fam\msbfam #1}}

\def \NN {\Bbb N}

\def \RR {\Bbb R}
\def \ZZ {\Bbb Z}

  \def\rightheadline{{\hfil{\ff
On some problems involving Hardy's function}\hfil\tenrm\folio}}

  \def\leftheadline{{\tenrm\folio\hfil{\ff
 Aleksandar Ivi\'c }\hfil}}
  \def\emptyheadline{\hfil}
  \headline{\ifnum\pageno=1 \emptyheadline\else
  \ifodd\pageno \rightheadline \else \leftheadline\fi\fi}

\topmatter
\title
On some problems involving Hardy's function
\endtitle
\author
 Aleksandar Ivi\'c
\endauthor
\address
Katedra Matematike RGF-a,
Universitet u Beogradu,  \DJ u\v sina 7,
11000 Beograd, Serbia.\bigskip
\endaddress
\subjclass 11 M 06
\endsubjclass
\email
{\tt
ivic\@rgf.bg.ac.rs,\enskip aivic\@matf.bg.ac.rs}
\endemail
\abstract
Some problems involving the classical Hardy function
$$
Z(t) := \zt\bigl(\chi(\hf+it)\bigr)^{-1/2}, \quad \z(s) = \chi(s)\z(1-s)
$$
are discussed. In particular we discuss the odd moments of $Z(t)$ and the
distribution of its positive and negative values.
\endabstract
\endtopmatter
\document
\head
1. Definition of Hardy's function
\endhead

\medskip
Let $Z(t)$ be the classical Hardy function (see e.g., [7]), defined as
$$
Z(t) := \zt\bigl(\chi(\hf+it)\bigr)^{-1/2},
$$
where $\chi(s)$ comes from the familiar functional equation for $\z(s)$ (see e.g., [7, Chapter 1]),
namely $\z(s) = \chi(s)\z(1-s)$, so that
$$
\chi(s) = 2^s\pi^{s-1}\sin(\hf \pi s)\G(1-s),\quad
\chi(s)\chi(1-s)=1.
$$
It follows  that
$$
\overline{\chi(\hf + it)} = \chi(\hf-it)= \chi^{-1}(\hf+it),
$$
so that $Z(t)\in\RR$ when $t\in\RR$
and $|Z(t)| =|\zt|$. Thus the zeros of $\z(s)$ on the ``critical line'' $\R s =1/2$
correspond to the real zeros of $Z(t)$, which makes $Z(t)$ an invaluable tool
in the study of the zeros of the zeta-function on the critical line.
Alternatively, if we use the symmetric form of the functional equation
for $\z(s)$, namely
$$
\pi^{-s/2}\z(s)\G(\hf s) = \pi^{-(1-s)/2}\z(1-s)\G(\hf(1-s)),
$$
then for $t\in \RR$ we obtain
$$
Z(t) = {\roman e}^{i\t(t)}\zt,\quad {\roman e}^{i\t(t)} := \pi^{-it/2}\frac
{\G({1\over4}+\hf it)}{|\G({1\over4}+\hf it)|},\quad \t(t) \in\RR.
$$
Hardy's original application of $Z(t)$ was to show that $\z(s)$ has infinitely
many zeros on the critical line $\R s =1/2$  (see e.g., E.C. Titchmarsh [23]).
The argument is briefly as follows.
Suppose on the contrary that, for $T\ge T_0$, the function $Z(t)$ does not change sign. Then
$$
\int_T^{2T}|Z(t)|\d t \;=\; \Bigl|\int_T^{2T}Z(t)\d t\Bigr|.\eqno(1.1)
$$
But it is not difficult to show (see [7, Chapter 9] and (2.4) for an even slightly
sharper bound) that
$$
\int_T^{2T}|Z(t)|\d t \;=\;\int_T^{2T}|\zt|\d t \;\gg\;T.\eqno(1.2)
$$
On the other hand, to bound the integral on the right-hand side of (1.1)
we can use the approximate functional equation (this is weakened
form of the so-called Riemann--Siegel formula; for a proof see [7] or [23])
$$
Z(t) = 2\sum_{n\le\sqrt{t/(2\pi)}}n^{-1/2}\cos\Biggl(t\log{\sqrt{t/(2\pi)}\over n}
- {t\over2} - {\pi\over8}\Biggr) + O\Bigl(\frac{1}{t^{1/4}}\Bigr).
$$
If this expression is integrated and the second derivative test is applied (see [7] or [23])
it follows that
$$
\int_T^{2T}Z(t)\d t \;=\; O(T^{3/4}).\eqno(1.3)
$$
Thus from (1.1)--(1.3) we obtain
$$
T \;\ll\; \int_T^{2T}|Z(t)|\d t \;\ll\;T^{3/4},
$$
which is a contradiction. This proves that $\z(s)$ has infinitely
many zeros on the critical line. Later Hardy refined his argument to show that $\z(s)$
has $\gg T$ zeros $\b + i\gamma$ satisfying $0 < \gamma \le T$. A. Selberg (see [22] or [23]
for a proof) improved this bound to $\gg T\log T$, which is one of the most important results
of analytic number theory. His method was later used and refined by many mathematicians.
The latest result is by S. Feng [1], who proved that at least 41.73\% of the zeros of $\z(s)$
are on the critical line and at least 40.75\% of  the zeros of $\z(s)$
are simple and on the critical line.

\medskip
\head
2. The distribution of values of Hardy's function
\endhead
\medskip

One of the  aims of this paper is to study the distribution of positive and negative
values of $Z(t)$. To this end let
$$
I_+(T) := \int_{T,Z(t)>0}^{2T}Z(t)\d t, \quad I_-(T) := \int_{T,Z(t)<0}^{2T}Z(t)\d t
$$
and
$$\eqalign{
&{\Cal J}_+(T) := \mu\Bigl\{T< t\le 2T\,:\, Z(t)>0\Bigr\},\cr&
{\Cal J}_-(T) := \mu\Bigl\{T< t\le 2T\,:\, Z(t)<0\Bigr\},\cr}\eqno(2.1)
$$
where $\mu(\cdot)$ denotes measure. Note that by the author's result [10] we have
$$
I_+(T) + I_-(T) = \int_T^{2T}Z(t)\d t = O_\e(T^{1/4+\e}),\eqno(2.2)
$$
which significantly improves on Hardy's bound $O(T^{3/4})$ (cf. (1.3)).
Here and later ``$\e$''  denotes arbitrarily small constants, not necessarily
the same ones at each occurrence, while $F = O_\e(G)$ means that the $O$-constant
depends on ``$\e$''.
Later M. Jutila [13] and M.A. Korolev [17] independently removed the ``$\e$'' in (2.2)
and proved that the integral is actually $\Omega_\pm(T^{1/4})$.
As usual,
$f = \Omega_+(g)$ means that $\limsup f/g > 0$, $f = \Omega_-(g)$ means that $\limsup f/g <0,$
and $f = \Omega_\pm (g)$ that both $f = \Omega_+(g)$ and $f = \Omega_-(g)$ hold.
Thus the problem of the true order of the integral of Hardy's function is settled,
up to the numerical values of the constants that are involved.
The primitive of Hardy's function will be discussed in more detail in the next section.

\medskip
Returning to the discussion on $I_\pm(T)$, note that we have
$$
\eqalign{
I_+(T) - I_-(T) &=  \int_{T,Z(t)>0}^{2T}Z(t)\d t - \int_{T,Z(t)<0}^{2T}Z(t)\d t\cr&
= \int_T^{2T}|Z(t)|\d t =  \int_T^{2T}|\zt|\d t,\cr} \eqno(2.3)
$$
since $|Z(t)| = |\zt|$. But from K. Ramachandra [21] we know that
$$
T(\log T)^{1/4} \;\ll\; \int_T^{2T}|\zt|\d t \;\ll\; T(\log T)^{1/4},\eqno(2.4)
$$
and from (2.2)--(2.3) it follows that
$$
I_+(T) = \hf \int_T^{2T}|\zt|\d t + O_\e(T^{1/4+\e}),
$$
$$
-I_-(T) = \hf \int_T^{2T}|\zt|\d t + O_\e(T^{1/4+\e}),
$$
and as mentioned, by [13] or [17] we can dispense with ``$\e$'' in the above two formulas.
In view of (2.4) we obtain then

\medskip
THEOREM 1. {\it We have}
$$\eqalign{&
T(\log T)^{1/4} \;\ll\; I_+(T)\ll T(\log T)^{1/4},\cr&
 T(\log T)^{1/4}\; \ll\; -I_-(T)\ll T(\log T)^{1/4}.\cr}
\eqno(2.5)
$$

\medskip\no
If one could sharpen (2.4) to an asymptotic formula, then we could sharpen (2.5)
and solve the following

\smallskip
{\bf Problem 1}. Prove that there is a constant $A>0$ such that
$$
I_+(T) = (A+ o(1))T(\log T)^{1/4},  \quad -I_-(T) = (A+ o(1))T(\log T)^{1/4}\quad(T\to\infty).
$$

\medskip
\head
3. The primitive of Hardy's function
\endhead
\medskip

Let us define
$$
F(T) := \int_0^T Z(t)\d t, \eqno(3.1)
$$
so that $F(t)$ is a primitive function of Hardy's function $Z(t)$. In
the aforementioned work [13] M. Jutila actually proved the following result.

\medskip
THEOREM A. {\it Let $T$ be a large positive number and write $\sqrt{T/(2\pi)} =
L + \t$ with $L\in\NN$ and $0\le \t < 1$. Define
$$
\t_0 \;=\; \min\Bigl(|\t - 1/4|,\,|\t - 3/4|\Bigr).
$$
Then for $\t_0\ne0$ we have
$$
F(T) = \Bigl({T\over2\pi}\Bigr)^{1/4}(-1)^LK(\t) + O(T^{1/6}\log T)
+ O\left(\min(T^{1/4},\,T^{1/8}\t_0^{-3/4})\right),
\eqno(3.2)
$$
where $K(x) = 0$ for $0\le x < 1/4$ and $ 3/4<x\le1$, $K(x) = 2\pi$ for $1/4<x<3/4$,
and further}
$$
\eqalign{
F(T)& = \Bigl({T\over2\pi}\Bigr)^{1/4}(-1)^L\,{4\pi\over3} + O(T^{1/6}\log T)
\quad{\roman {for}}\,\t=1/4,\cr
F(T)&= \Bigl({T\over2\pi}\Bigr)^{1/4}(-1)^L\,{2\pi\over3} + O(T^{1/6}\log T)
\quad{\roman {for}}\,\t=3/4.\cr}
$$
The main contribution to $F(T)$ comes from
the expression
 $$
 F_1(T) := \Bigl({T\over2\pi}\Bigr)^{1/4}(-1)^LK(\t)
 = 2\pi\Bigl({T\over2\pi}\Bigr)^{1/4}(-1)^{[\sqrt{T/(2\pi)}]}
  $$
 when $1/4 < \t = \{\sqrt{T/(2\pi)}\} < 3/4$, and $F_1(T)=0$ otherwise.
 In [14] Jutila derived another expression for $F(T)$, different from (3.1),
  where $F_1(T)$ is
 replaced by a smoothed expression involving the integral of $K(\t+t)-K(t)$, plus
 an error term which is only $O(T^{1/6}\log T)$.
 On  the average $F_1(T)$ is small, as shown by the following

 \medskip
 THEOREM 2. {\it We have}
 $$
 \eqalign{
 &\int_0^{T}F_1(t)\d t \;=\;O(T^{3/4}),\cr&
 \int_0^{T}F_1(t)\d t \;=\;\Omega_\pm(T^{3/4}).\cr}
\eqno(3.3)
 $$

\medskip
{\bf Proof.} Setting  $u(x) = (-1)^{[x]}K(\{x\})$ ($[x]$ is the integer part of $x$
and $\{x\}$ is its fractional part) we see that this is an odd function
satisfying $u(x+2)  = u(x)$, and it follows that
we have the Fourier series (for $x- \txt{1\over4}\not\in\ZZ,\;x- \txt{3\over4}\not\in\ZZ)$
$$
u(x) = 4\sum_{n=1}^\infty  {\cos({\pi n\over4})-\cos({3\pi n\over4})\over n}\sin(\pi nx).
$$
Setting $a(n) := \cos({\pi n\over4})-\cos({3\pi n\over4})$ we see that $a(2k)=0$
and that $a(n) = \sqrt{2}$ for $n \equiv 1,7\,$(mod 8) and $a(n) = -\sqrt{2}$
for $n \equiv 3,5\,$(mod 8). This gives
$$
u(x) = 4\sum_{k=1}^\infty {a(2k-1)\over2k-1}\sin\bigl(\pi(2k-1)x\bigr)
\quad(x- \txt{1\over4}\not\in\ZZ,\;x- \txt{3\over4}\not\in\ZZ).
$$
Hence with the change of variable $t =2\pi x^2$ and integration by parts
we obtain, since the above series
is boundedly convergent and can be integrated termwise,
$$
\eqalign{&
\int_0^TF_1(t)\d t = \int_0^T\left({t\over2\pi}\right)^{1/4}u
\left(\bigl({t\over2\pi}\bigr)^{1/2}\right)\d t\cr&
= 16\pi\int_0^{\sqrt{T/2\pi}}x^{3/2}\sum_{k=1}^\infty {a(2k-1)\over2k-1}
\sin\bigl(\pi(2k-1)x\bigr)\d x\cr&
= 16\pi\sum_{k=1}^\infty {a(2k-1)\over2k-1}\int_0^{\sqrt{T/2\pi}}x^{3/2}
\sin\bigl(\pi(2k-1)x\bigr)\d x\cr&
= -16\left({T\over2\pi}\right)^{3/4}\sum_{k=1}^\infty {a(2k-1)\over(2k-1)^2}
\cos\left(\pi(2k-1)\sqrt{{T\over2\pi}}\right) + O(T^{1/4}).
\cr}
$$
Since the last series is absolutely convergent
 we obtain the $O$-estimate of (3.3). The Omega-results
follow if we take $T = 2\pi({3\over 4}+2m)^2$ and $T = 2\pi({1\over4}+2m)^2$ with $m\in\NN$.
In fact the last $O$-term
stands for a function that is $O(T^{1/4})$ and also $\Omega_\pm(T^{1/4})$.
This proves (3.3). However the true order of the primitive of $F(T)$ remains
elusive, since it is not obvious how small will be, when integrated,
 the expression standing for the $O$-term in Jutila's expression for $F(T)$.

\medskip
{\bf Problem 2}. What is the true order of
$$
\int_0^{T}F(t)\d t?
$$
Is it perhaps true that
$$
\int_0^{T}F(t)\d t\;=\; O(T^{3/4}),\quad\int_0^{T}F(t)\d t\;=\; \Omega_\pm(T^{3/4})?
$$

\medskip
\head
4. The cubic moment of Hardy's function
\endhead
\medskip

The analogous problem when $Z(t)$ in $I_\pm(T)$ is replaced by $Z^3(t)$ is much harder.
Namely there is an old problem of mine (posed in Oberwolfach, 2003), which is stated here
as

\smallskip
{\bf Problem 3}. Does there exist  a constant $0<c<1$ such that
$$
\int_0^TZ^3(t)\d t \;=\;O(T^c)?\eqno(4.1)
$$

\smallskip

\noindent
One may naturally ask for bounds for higher moments of $Z(t)$. However, only odd
moments are interesting in this context, because of their oscillating nature.
Unfortunately when $k>4$
not much, in general, is known on the moments of $|\zt|^{k}$  (see e.g., [7]).
I have not been able to prove (4.1) yet, although I feel that there must be a lot
of cancelation between the positive and negative values of $Z(t)$ and I am
certain that it must be true. What can be proved is (see [12])
$$
\int\limits_T^{2T}Z^3(t)\d t= 2\pi\sqrt{2\over 3}
\sum_{({T\over 2\pi})^{3/2}\le n\le ({T\over\pi})^{3/2}}
d_3(n)n^{-{1\over6}}\cos\bigl(3\pi n^{2\over 3}+{\txt{1\over8}}\pi\bigr)
+O_\e(T^{3/4+\e}),\eqno(4.2)
$$
where as usual the divisor function $d_3(n)$ (generated by $\z^3(s)$) denotes the
number of ways $n$ can be written as a product of three factors. The difficulty is in
the estimation of the exponential sum on the right-hand side of (4.2). I can show that
$$
\sum_{N\le n\le 2N}
d_3(n)n^{-1/6}\cos\bigl(3\pi n^{2\over 3}+{\txt{1\over8}}\pi\bigr) \ll_\e N^{2/3+\e}.\eqno(4.3)
$$
However, this just gives the bound $O_\e(T^{1+\e})$ for the integral on the left-hand
side of (4.2). This is unfortunately weak, since by the Cauchy-Schwarz inequality for
integrals we easily obtain a better result, namely
$$
\Bigl|\,\int\limits_T^{2T}Z^3(t)\d t\,\Bigr| \le \left(\int\limits_T^{2T}|\zt|^2\d t
\int\limits_T^{2T}|\zt|^4\d t\right)^{1/2} \ll T(\log T)^{5/2}\eqno(4.4)
$$
on using the well-known elementary bounds (see e.g., [7] or [23])
$$
\int_T^{2T}|\zt|^2\d t \;\ll\; T\log T,\quad \int_T^{2T}|\zt|^4\d t \;\ll\; T\log^4T.
$$
For the cubic moment of $Z(t)$ I know of no better bound than (4.4),
and any improvement thereof would be interesting.
A nice feature of the exponential sum in (4.3) is that it is ``pure'' in the sense that
the argument of the cosine depends only on $n$, and not on any other quantity.
A reasonable conjecture is that
$$
\int_0^TZ^3(t)\d t \;=\;O_\e(T^{3/4+\e}), \quad \int_0^TZ^3(t)\d t = \Omega_\pm(T^{3/4}).
$$
Actually (4.3) holds if $d_3(n)$ is replaced by $d_k(n)$ (the number of ways $n$ can be
written as a product of $k$ factors, so that $d_1(n)\equiv1$), provided that we have
$$
\int_0^T|\zt|^k\d t \;\ll_{\e,k}\; T^{1+\e}\qquad(k\in\NN).\eqno(4.5)
$$
In (4.5) it is assumed that $k$ is fixed. If (4.5) holds for any $k\;(\ge k_0)$, then this
is equivalent to the Lindel\"of hypothesis that $|\zt| \ll_\e |t|^\e$ (see e.g., E.C.
Titchmarsh [23]).

\smallskip
The assertion related to (4.3) is contained in

\medskip
THEOREM 3. {\it If} (4.5) {\it holds for some fixed $k\in\NN$,  then}
$$
\sum_{N\le n\le 2N}
d_k(n)n^{-1/6}\cos\bigl(3\pi n^{2\over 3}+{\txt{1\over8}}\pi\bigr)\; \ll_\e\; N^{2/3+\e}.
\eqno(4.6)
$$

\medskip

\no The bound (4.5) is at present known to hold only when $k\le 4$, in which case
we have a non-trivial result. It is unclear whether (4.6), in the general case, has
an arithmetic meaning as in the case $k=3$ (cf. (4.2)).

\medskip
{\bf Proof of Theorem 3.}  From the Perron inversion formula (see e.g., [7, Appendix]) we have
$$
\sum_{n\le x} d_k(n)= {1\over2\pi i}\int_{1+\e-iT}^{1+\e+iT}{\z^k(w)x^w\over w}\d w
+ O_\e(N^\e)\qquad(T\asymp x, x\asymp N).\eqno(4.7)
$$
We replace the segment of integration by $[\hf-iT, \hf+iT]$, passing over the pole
$w=1$ of order $k$ of $\z^k(w)$. The residue is $xP_{k-1}(\log x)$, with $P_{k-1}(\log x)$
a polynomial of degree $k-1$ in $\log x$. It follows by the residue theorem that
$$
\D_k(x) := \sum_{n\le x} d_k(n) - xP_{k-1}(\log x)
= {1\over2\pi i}\int_{{1\over2}-iT}^{{1\over2}+iT}{\z^k(w)x^w\over w}\d w
+ R_k(x,T),\eqno(4.8)
$$
say. Here $\D_k(x)$ is the error term in the asymptotic formula for the summatory
function of $d_k(n)$, and $R_k(x,T)$ stands for the error term in (4.7) plus
the contribution over the segments $[\hf\pm iT,\,1+\e\pm iT]$, so that
$$
R_k(x,T) \ll_\e \int_{1\over2}^{1+\e}|\z(\s + iT)|^kx^\s T^{-1}\d \s + N^\e.\eqno(4.9)
$$
If we integrate (4.9) over $T$ from $T_1$ to $2T_1$, similarly as was done in [11], we
shall find a value of $T \;(\asymp N)$ such that $R_k(x,T) \ll_{\e,k} N^\e$. It is
actually this value of $T$ that is taken in (4.7). In this
process we need the bound, which follows by (4.5) and the convexity of mean
values (see [7, Chapter 8]). This is
$$
\int_0^T|\z(\s+it)|^k\d t \;\ll_{\e,k}\; T^{1+\e}\qquad(\hf \le \s \le 1).
$$
We write now the sum in (4.6) as
$$
\int_N^{2N}x^{-1/6}\cos\bigl(3\pi x^{2\over 3}+{\txt{1\over8}}\pi\bigr)
\d\Bigl(\sum_{n\le x} d_k(n) \Bigr)
$$
and use (4.8). The contribution of $xP_{k-1}(\log x)$, by the first derivative test
([7, Lemma 2.1]), will be $\ll N^{1/6+\e}$. In the portion pertaining to $\D_k(x)$
we integrate by parts the term $R_k(x,T)$. Since $R_k(x,T) \ll_{\e,k} N^\e$, trivial
estimation will yield a contribution which is $\ll_\e N^{1/2+\e}$,
and this is probably optimal.
After differentiating the remaining expression over $x$ and writing $w = \hf + iv$
and the cosine as a
sum of exponentials, we are left with a multiple of
$$
\int_{-T}^T \z^k(\hf + iv)\left\{\int_N^{2N}x^{-2/3}\exp\left(iv\log x
\pm 3\pi ix^{2/3}\right)\d x\right\}\d v.
$$
The function in the exponential is $iF_\pm(x)$ with
$$
F_\pm(x) = F_\pm(x,v) := v\log x \pm 3\pi x^{2/3},
$$
so that
$$
F_\pm'(x) = {v\over x} \pm 2\pi x^{-1/3},\quad F_\pm''(x) = -{v\over x^2}
\mp {2\over3}\pi x^{-1/3}.
$$
Suppose $v>0$, the other case being analogous. The saddle point $x_0$ (root
of $F_+'(x)=0$ in this case) is $x_0 = (v/(2\pi))^{3/2}$, and $x_0\in[N,2N]$
for $v \asymp N^{2/3}$. When this is not satisfied, the contribution is,
by the first derivative test, $\ll_\e N^{1/3+\e}$. For $v \asymp N^{2/3}$ we have
$$
|F_+''(x_0)|^{-1/2} \;\asymp\; v \;\asymp\; N^{2/3}.
$$
Hence by the second derivative test and (4.5) we have
that this contribution is,
with suitable constants $0 <C_1<C_2$,
$$
\ll \int_{C_1N^{2/3}}^{C_2N^{2/3}}|\z(\hf+iv)|^k\d v\cdot N^{-2/3}N^{2/3} \ll_\e
N^{2/3+\e}.
$$
Therefore (4.6) follows, and Theorem 3 is proved. If one wants to have explicitly
the main term produced by this method, one can use the saddle point method
(see e.g., [7, Chapter 2]).

\medskip
\head
5. Further problems on the distribution of values
\endhead
\medskip

If (4.1) holds with some $0<c<1$ then, similarly as in the case of $I_\pm(T)$, we obtain
$$
\eqalign{&
\int_{T,Z(t)>0}^{2T}Z^3(t)\d t \;=\; \hf \int_T^{2T}|\zt|^3\d t + O(T^c),\cr&
\int_{T,Z(t)<0}^{2T}Z^3(t)\d t \;=\; \hf \int_T^{2T}|\zt|^3\d t + O(T^c).\cr}
\eqno(5.1)
$$
Note that no asymptotic formula exists yet for the integral on the right-hand side
of (5.1). In general it is conjectured that (known to be true only when $k=2,4$)
$$
\int_0^T|\zt|^k\d t \;=\;(c_k + o(1))T(\log T)^{k^2/4}\qquad(k\in\NN,\,\,T\to\infty),
$$
where (see  J. P. Keating and N.C. Snaith [16]), with $d_k(n)$ generated by $\z^k(s)$,
the value of $c_k$  for even $k$ is
$$
c_k = {a_kg_k\over\G(1+k^2)},\quad a_k = \prod_{p}\Bigl(1-\frac{1}{p}\Bigr)^{k^2}
\sum_{j=0}^\infty d^2_k(p^j),\quad g_k = (k^2)!\prod_{j=0}^{k-1}{j!\over(j+k)!},
$$
and the product is taken over all primes $p$.
In the case $k=3$ all that is currently known is
$$
T(\log T)^{9/4} \ll \int_0^T|\zt|^3\d t \ll T(\log T)^{5/2}.
$$
The upper bound follows as in (4.4), and the lower bound is a
special case of the general lower bound (see e.g., K. Ramachandra [21])
$$
\int_0^T|\zt|^k\d t \;\gg_k\; T(\log T)^{k^2/4}\qquad(k\in\NN).
$$

\medskip

The problems of the asymptotic evaluation of $I_\pm(T)$ are clearly connected to the distribution
of the values of $\zt$ and $|\zt|$. In [9] I proved, using results of A. Selberg (see
D.A. Hejhal [6]) on the distribution of values of $\z(s)$, that
$$
\mu(A_c(T)) = {T\over2} +
O\left(T{(\log\log \log T)^2\over\sqrt{\log\log T}}\right),\eqno(5.2)
$$
where $c>0$ is any constant and
$$
A_c(T) := \Bigl\{ 0 < t \le T : |\zt| \le c \Bigr\}.
$$
Unfortunately one does not see how to put (5.2) to use in connection with the
distribution of positive values of $Z(t)$. It is hard to determine asymptotically
the order of ${\Cal J}_+(T), {\Cal J}_-(T)$, defined by (2.1). We formulate
the following

\smallskip
{\bf Problem 4}. Is it true that there exist constants $A_+>0, A_->0$ such that
$$
{\Cal J}_+(T) := \mu\Bigl\{T< t\le 2T\,:\, Z(t)>0\Bigr\} = \bigl(A_+ + o(1)\bigr)T\quad(T\to\infty),
\eqno(5.3)
$$
$$
{\Cal J}_-(T) := \mu\Bigl\{T< t\le 2T\,:\, Z(t)<0\Bigr\} = \bigl(A_- + o(1)\bigr)T\quad(T\to\infty)?
\eqno(5.4)
$$

\smallskip

Obviously $A_+ + A_- = 1$ (if $A_+, A_-$ exist). The asymptotic formula (5.2) gives
rise to the thought that in that case maybe $A_+ = A_- = 1/2$. On the other hand,
things may not be that simple. If one assumes the famous Riemann Hypothesis (that
all complex zeros of $\z(s)$ lie on $\R s = 1/2$) and the simplicity of zeta zeros
(these very strong conjectures seem to be independent in the sense that it is not
known whether either of them implies the other one) then (since $Z(0) = -1/2$)
$$
\mu\Bigl\{ T < t \le 2T : Z(t)>0\Bigr\} = \sum_{T<\gamma_{2n}\le 2T}(\gamma_{2n}
- \gamma_{2n-1}) + O(1),\eqno(5.5)
$$
where $0 < \gamma_1 < \gamma_2 < \ldots$ are the ordinates of complex zeros of $\z(s)$.
The  sum in (5.5) is connected to the sum ($\a\ge0$ is fixed)
$$
\sum\nolimits_\a(T) : = \sum_{\gamma_n\le T}(\gamma_n-\gamma_{n-1})^\a,
$$
which was investigated in [8]. The sum $\sum_\a(T)$ in turn  can be connected to the
Gaussian Unitary Ensemble hypothesis (see A.M. Odlyzko [19], [20])
and the pair correlation conjecture of H.L. Montgomery [18]. Both of these conjectures
assume the Riemann Hypothesis and e.g., the former states that, for
$$
0 \le \a < \b < \infty, \quad \delta_n = \frac{1}{2\pi}(\gamma_{n+1}-\gamma_n)\log\left(
\frac{\gamma_n}{2\pi}\right),
$$
we have
$$
\sum_{\gamma_n\le T,\delta_n\in[\a,\b]}1 = \left(\int_\a^\b p(0,u)\d u + o(1)\right)
{T\over2\pi}\log\left({T\over2\pi}\right)\quad(T\to\infty),
$$
where $p(0,u)$ is a certain probabilistic density, given by complicated functions
defined in terms of prolate spheroidal functions. In fact, in [8] I have proved that,
if the RH and the Gaussian Unitary Ensemble hypothesis hold, then for $\a\ge0$ fixed
and $T\to\infty$,
$$
\sum\nolimits_\a(T) = \Bigl(\int_0^\infty p(0,u)u^\a\d u + o(1)\Bigr)
{\left({2\pi\over \log\bigl({T\over2\pi}\bigr)-1}\right)}^{\a-1}T.
$$
\smallskip
Also note that, since $\R\log\zt
= \log|Z(t)|$, a classical result of A. Selberg (see [22]) gives, for any real $\a<\b$,
$$
\lim_{T\to\infty}{1\over T}\mu\Biggl\{\,t : t\in [T,2T],\,
\a < {\log|Z(t)|\over\sqrt{\hf\log\log T}} < \b\,\Biggr\}
= {1\over\sqrt{2\pi}}\int_\a^\b {\roman e}^{-{1\over2}x^2}\d x,
$$
but here we are interested in the distribution of values of $Z(t)$ and not $|Z(t)|$.
\smallskip
Recently J. Kalpokas and J. Steuding [15] proved, among other things, that for \break
$\phi\in [0,\pi)$,
$$
\sum_{0<t\le T,\zt\in {\roman e}^{i\phi}\RR}\zt = \Bigl(2{\roman e}^{i\phi}
\cos\phi\Bigr) {T\over2\pi}\log{T\over2\pi{\roman e}} + O_\e(T^{1/2+\e}),\eqno(5.6)
$$
and an analogous result holds for the sums of $|\zt|^2$. It is unclear whether (5.6)
and the other approaches mentioned above can be put to
use in connection with our problems.

\smallskip
Although (5.3) and (5.4) seem difficult to prove, one can at least show that
$$
{\Cal J}_+(T) \;\gg\; T(\log T)^{-1/2},\eqno(5.7)
$$
and a similar bound for ${\Cal J}_-(T)$. Namely from (2.5) we have, by
the Cauchy-Schwarz inequality,
$$
T(\log T)^{1/4} \ll I_+(T) \le \left(\int_{T,Z(t)>0}^{2T}1\d t\,
\int_T^{2T}|\zt|^2\d t\right)^{1/2},
$$
which easily gives (5.7), since
$$
\int_{T,Z(t)>0}^{2T}1\d t \;=\; \mu\Bigl\{ T < t \le 2T : Z(t)>0\Bigr\} = {\Cal J}_+(T).
$$

\medskip
\head
6. The distribution of values of $E(T)$
\endhead
\medskip

The problem analogous to the evaluation of $I_\pm(T)$ can be considered for
$E(T)$ (see [7, Chapter 15]), the error term in the mean square formula for $|\zt|^2$. Namely
one has (here $C_0= - \G'(1)$ is Euler's constant) the defining relation
$$
E(T) = \int_0^T|\zt|^2\d t - T\log{T\over2\pi} -(2C_0-1)T.
$$
Thus the function $E(T)$ is continuous for $T>0$,
and it is known that $E(T) = O(T^\a)$ with $1/4 \le \a < 1/3$ (see [7, Chapter 15]).
Moreover its mean value is $\pi$, as  J.L. Hafner and I showed in [2]-[3] that if we define
$$
G(T) \;:=\; \int_0^T(E(t)-\pi)\d t,
$$
one can then obtain precise expressions for $G(T)$ which, in particular, yield
$$
 G(T) \;=\; O(T^{3/4}),\quad G(T) \;=\; \Omega_\pm(T^{3/4}).
$$
Let us also define
$$
J_+(T) \;:=\; \int_{T,E(t)>\pi}^{2T}(E(t)-\pi)\d t,
\quad
J_-(T) \;:=\; \int_{T,E(t)<\pi}^{2T}(E(t)-\pi)\d t.
$$
Then, similarly to the discussion on $I_\pm(T)$, we obtain
$$
J_+(T) = \hf\int_T^{2T}|E(t)-\pi|\d t + O(T^{3/4}),\quad
-J_-(T) = \hf\int_T^{2T}|E(t)-\pi|\d t + O(T^{3/4}).\eqno(6.1)
$$
Note that D.R. Heath-Brown [4] proved that the moments
$$
D_k \;:=\; \lim_{X\to\infty}X^{-1-k/4}\int_0^X|E(t)|^k\d t\eqno(6.2)
$$
exist for any real $k\in[0,\,9]$, and so do the odd moments of $E(t)$ for $k=1,3,5,7$ or 9.
If $D_k$ exists for some $k\ge1$, then $D_k > 0$. Namely from a result of
D.R. Heath-Brown and K.-M. Tsang [5] it follows that
$$
\int_0^X|E(t)|\d t  \;\gg\; X^{5/4},\eqno(6.3)
$$
so that $D_1 >0$. If $k>1$, then by H\"older's inequality
$$
\int_0^X|E(t)|\d t \;\ll_k\; {\left(\int_0^X|E(t)|^k\d t\right)}^{1/k}X^{1-1/k},
$$
so that (6.3) yields
$$
\int_0^X|E(t)|^k\d t  \;\gg_k\; X^{1+k/4},
$$
implying that $D_k>0$ if $D_k$ exists. Since
$$
\int_T^{2T}|E(t)|\d t - \pi T\le \int_T^{2T}|E(t)-\pi|\d t\le \int_T^{2T}|E(t)|\d t +\pi T,
\eqno(6.4)
$$
then using (6.2) (for $k=1)$ once with $X=T$ and
once with $X = 2T$, we obtain from (6.1) and (6.4)

\medskip
THEOREM 4. {\it As $T\to\infty$ we have, for some constant} $C>0$,
$$
\eqalign{
J_+(T) &\;=\;(C+o(1))T^{5/4},\cr -J_-(T) &\;=\;(C+o(1))T^{5/4}.\cr}
$$

\medskip
We conclude by stating two related problems.
\medskip
{\bf Problem 5}. Does there exist a constant $B>0$ such that
$$
\int_0^T |E(t)|\d t \;=\; BT^{5/4} + O(T)?
$$

\smallskip
{\bf Problem 6}. Do there exist a constants $B_+, B_->0$ such that, for $T\to\infty$,
$$\eqalign{&
\mu\Bigl\{ \,0 \le t \le T \,:\, E(t) > \pi\,\Bigr\} = \bigl(B_+ +o(1)\bigr)T,\cr&
\mu\Bigl\{ \,0 \le t \le T \,:\, E(t) < \pi\,\Bigr\} = \bigl(B_- +o(1)\bigr)T?\cr}
$$
Obviously $B_+ + B_- = 1$ (if $B_+, B_-$ exist), and  maybe $B_+ = B_- = 1/2$.
This is analogous to (5.3)--(5.4). Generalizations to the distribution of positive and
negative values of other functions of arithmetic interest may be clearly considered as well.
\vskip2cm
\Refs
\smallskip

\item{[1]}  Feng S., Zeros of the Riemann zeta-function
on the critical line, preprint available at \break arXiv:1003.0059

\item{[2]}  Hafner J.L. and  Ivi\'c A., On some mean value results for the
  Riemann zeta-function,  Proceedings International Number Theory Conference
  Qu\'ebec 1987, Walter de Gruyter and Co., 1989, Berlin - New York, 348-358

\item{[3]} Hafner J.L. and  Ivi\'c A., On the mean square of the Riemann zeta-function
on the critical line,  J. Number Theory   {\bf 32}(1989), 151-191

 \item{[4]}  Heath-Brown D.R., The distribution and moments of the error term
 in the Dirichlet divisor problems,  Acta Arith. {\bf60}(1992), 389-415

 \item{[5]} Heath-Brown D.R. and Tsang K., Sign changes of $E(T), \D(x)$ and $P(x)$,
  J. Number Theory {\bf49}(1994), 73-83

\item{[6]}  Hejhal D.A.,  On a result of Selberg concerning
zeros of linear combinations of L-functions,
 Int. Math. Res. Not. 2000, No. 11, 551-577 (2000)

\item{[7]} Ivi\'c A., The Riemann zeta-function, John Wiley \&
Sons, New York 1985 (2nd edition.  Dover, Mineola, New York, 2003)

\item {[8]} Ivi\'c A., On sums of gaps between the zeros of
$\zeta(s)$ on the critical line,
 Univ. Beograd. Publ. Elektrotehn. Fak. Ser. Mat. {\bf6}(1995), 55-62

\item {[9]} Ivi\'c A., On small values of the Riemann zeta-function on the critical line
and gaps between zeros, Lietuvos Mat. Rinkinys {\bf42}(2002), 31-45

\item {[10]} Ivi\'c A., On the integral of Hardy's function,  Arch. Mathematik
{\bf83}(2004), 41-47

\item {[11]} Ivi\'c A., On the mean square of the divisor function in short intervals,
 Journal de Th\'eorie des Nombres de Bordeaux {\bf21}(2009), 195-205

\item {[12]} Ivi\'c A., On the Mellin transforms of powers of Hardy's function,
 Hardy-Ramanujan Journal {\bf33}(2010), 32-58

\item{[13]}  Jutila M., Atkinson's formula for Hardy's function,
J. Number Theory {\bf129}(2009), 2853-2878

\item{[14]} Jutila M., An asymptotic formula for the primitive of Hardy's
function,  Arkiv Mat. \break DOI:10.1007/s11512-010-0122-4

\item{[15]}  Kalpokas J. and  Steuding J., On the value distribution of the Riemann
zeta-function on the critical line, preprint available at arXiv:0907.1910

\item{[16]}  Keating J.P. and   Snaith N.C., Random Matrix Theory
and $L$-functions at $s=1/2$,   Comm. Math. Phys.{\bf214}(2000), 57-89

\item{[17]}  Korolev M.A., On the integral of Hardy's function $ Z(t)$,
 Izv. Math. {\bf72}, No. 3, 429-478 (2008); translation from  Izv.
Ross. Akad. Nauk, Ser. Mat. {\bf72}, No. 3, 19-68 (2008)

\item{[18]}  Montgomery H.L., The pair correlation of zeros of the
zeta-function,  Proc. Symp. Pure Math. {\bf24}, AMS, Providence 1973,
181-193

\item{[19]}  Odlyzko A.M., On the distribution  of spacings of zeros of
the zeta-function,  Math. Comp. {\bf48}(1987), 273-308

\item{[20]}  Odlyzko A.M., The $10^{20}$-th zero of the Riemann zeta-function
and 175 million of its neighbors, preprint available at

\item{} {\tt http://www.dtc.umn.edu/$\thicksim$odlyzko/unpublished/zeta.10to20.1992.pdf}

\item{[21]}  Ramachandra K., On the mean-value and omega-theorems
for the Riemann zeta-function,  Tata Inst. of Fundamental Research
(distr. by Springer Verlag, Berlin etc.), Bombay, 1995

\item{[22]} Selberg A., Selected papers, Vol. 1,  Springer Verlag, Berlin, 1989

\item{[23]}  Titchmarsh E.C., The theory of the Riemann
zeta-function (2nd edition),   University Press, Oxford, 1986

\endRefs

\enddocument

\bye